\documentclass[letterpaper, 10pt, conference]{ieeeconf}

	\IEEEoverridecommandlockouts                              
	\PassOptionsToPackage{noend}{algpseudocode}
	\PassOptionsToPackage{breakspaces,capitalize}{clevethm}
	\let\labelindent\relax

	\makeatletter
		\let\NAT@parse\undefined
	\makeatother
	
	\usepackage{amsthm}
	\usepackage[nospace]{cite}
	\usepackage[noload={caption,subfig,subcaption}]{myPreamble}
	\hypersetup{
		breaklinks,
		plainpages    = false,
		colorlinks    = true,
		urlcolor      = BlueIMT,
		linkcolor     = MidnightBlue,
		citecolor     = Green,
		hypertexnames = false,
	}

		\CreateNewTheorem[Fact]{fact}[dummythm]{Fact}
	
	\newcommand{\TheThanks}{%
		\textsuperscript1Andreas Themelis and Panagiotis Patrinos are with the \TheAddressKU.
		This work was supported by the \emph{Research Foundation Flanders (FWO)} research projects G086518N, G086318N, and G0A0920N;
		\emph{Research Council KU Leuven} C1 project No. C14/18/068;
		\emph{Fonds de la Recherche Scientifique --- FNRS and the Fonds Wetenschappelijk Onderzoek --- Vlaanderen} under EOS project no 30468160 (SeLMA).%
		\newline
		\textsuperscript2Ben Hermans is with the Department of Mechanical Engineering, KU Leuven, and DMMS lab, Flanders Make, Leuven, Belgium.
		His research benefits from KU Leuven-BOF PFV/10/002 Centre of Excellence: Optimization in Engineering (OPTEC), from project G0C4515N of the Research Foundation - Flanders (FWO - Flanders), from Flanders Make ICON: Avoidance of collisions and obstacles in narrow lanes, and from the KU Leuven Research project C14/15/067: B-spline based certificates of positivity with applications in engineering.%
		\newline
		{\tt
			\{\href{mailto:andreas.themelis@kuleuven.be}{andreas.themelis},%
			\href{mailto:ben.hermans2@kuleuven.be}{ben.hermans2},%
			\href{mailto:panos.patrinos@kuleuven.be}{panos.patrinos}\}%
			\href{andreas.themelis@kuleuven.be,mailto:ben.hermans2@kuleuven.be,panos.patrinos@kuleuven.be}{@kuleuven.be}%
		}%
	}%
	
	\crefname{section}{{\S}}{{\S}}
	\Crefname{section}{Section}{Sections}
	\crefformat{section}{#2{\S}#1#3}
	\Crefformat{section}{#2Section #1#3}
	\crefname{figure}{Figure}{Figures}

	\makeatletter
		\newcommand{\f}{f}
		\newcommand{\g}{g}
		\newcommand{\h}{h}
		\newcommand{\G}{G}
		\renewcommand{\H}{H}
		\newcommand{\genf}{f}
		\newcommand{\cost}{\varphi}
		\newcommand{\Cost}{\Phi}
		\newcommand{\envgamma}{\gamma}
		\newcommand{\envg}{\g}
		\newcommand{\envh}{\h}
		\newcommandx{\env}[3][1={\envh},2={\envg},3={\envgamma}]{\operatorname{env}_{#3}^{#2,#1}}
	\makeatother

	\newcommand{\FBE}{\cost_\gamma^{\text{\sc fb}}}

	\renewcommand{\trans}[1]{#1^{{\top}}}
	\DeclareMathOperator{\symm}{sym}
	
\showgrayfalse

\title{\LARGE\bf A new envelope function for nonsmooth DC optimization}
\author{%
	Andreas Themelis,\texorpdfstring{\textsuperscript1}{}
	Ben Hermans,\texorpdfstring{\textsuperscript2}{}
	and Panagiotis Patrinos\texorpdfstring{\textsuperscript1}{}%
	\texorpdfstring{\thanks{\TheThanks}}{}%
}

\renewcommand{\includetikz}[2][]{\includegraphics[#1]{Pics/Tikz/#2.pdf}}

\begin{document}

	\maketitle
	\begin{abstract}
	Difference-of-convex (DC) optimization problems are shown to be equivalent to the minimization of a Lipschitz-differentiable ``envelope''.
	A gradient method on this surrogate function yields a novel (sub)gradient-free proximal algorithm which is inherently parallelizable and can handle fully nonsmooth formulations.
	Newton-type methods such as L-BFGS are directly applicable with a classical linesearch.
	Our analysis reveals a deep kinship between the novel DC envelope and the forward-backward envelope, the former being a smooth and convexity-preserving nonlinear reparametrization of the latter.
	\end{abstract}

	\section{Introduction}
	We consider difference-of-convex (DC) problems
	\[\tag{P}\label{eq:P}
		\minimize_{s\in\R^p}\cost(s)\coloneqq\g(s)-\h(s),
	\]
	where \(\func{\g,\h}{\R^p}{\R\cup\set\infty}\) are proper, convex, lsc functions (with the convention \(\infty-\infty=\infty\)).
	DC problems cover a very broad spectrum of applications; a well detailed theoretical and algorithmic analysis is presented in \cite{tao1997convex},
	where the nowadays textbook algorithm DCA is presented that
	interleaves subgradient evaluations \(v\in\partial{\h}(u)\), \(u^+\in\partial\conj{\g}(v)\), aiming at finding a \emph{stationary} point \(u\), that is, a point satisfying
	\begin{equation}\label{eq:stationary}
		\partial\g(u)\cap\partial\h(u)\neq\emptyset,
	\end{equation}
	a relaxed version of the necessary condition \(\partial\h(u)\subseteq\partial\g(u)\) \cite{hiriarturruty1989convex}.
	As noted in \cite{an2017convergence}, proximal \emph{sub}gradient iterations are effective even in handling a nonsmooth nonconvex \(\g\) and a nonsmooth concave \(-\h\).
	Alternative approaches use the identity \(-\f(x)=\inf_y\set{\conj{\f}(y)-\innprod xy}\) involving the convex conjugate \(\conj{\f}\) to include an additional convex function \(\f\) as
	\begin{equation}\label{eq:P3}
		\minimize_{x\in\R^n}\g(x)-\h(x)-\f(x),
	\end{equation}
	and then recast the problem as
	\begin{equation}\label{eq:P3DC}
		\minimize_{x,y\in\R^n}{
			\Cost(x,y)
		{}\coloneqq{}
			\overbracket{\g(x)+\conj{\f}(y)}^{\G(x,y)}
			{}-{}
			\bigl(
				\overbracket{\vphantom{\conj{\f}}\h(x)+\innprod xy}^{\H(x,y)}
			\bigr)
		}.
	\end{equation}
	By adding and substracting suitably large quadratics, one can again obtain a decoupled DC formulation, showing that \eqref{eq:P} is in fact as general as \eqref{eq:P3}.
	When function \(\h\) is smooth (differentiable with Lipschitz gradient), a cornerstone algorithm for the ``convex\(+\)smooth'' formulation \eqref{eq:P3DC} is forward-backward splitting (FBS), amounting to gradient evaluations of the smooth component
	\(
		-\h(s)-\innprod st
	\)
	followed by proximal operations (possibly in parallel) on \(\g\) and \(\conj{\f}\).
	
	A detailed overview on DC algorithms is beyond the scope of this paper; the interested reader is referred to the exhaustive surveys in \cite{tao1997convex,horst1999dc,bacak2011difference} and references therein.
	Most related to our approach, \cite{banert2019general} analyzes a Gauss-Seidel-type FBS in the spirit of the PALM algorithm \cite{bolte2014proximal}, and \cite{liu2017further} exploits the interpretation of FBS as a gradient-type algorithm on the \emph{forward-backward envelope} (FBE) \cite{patrinos2013proximal,stella2017forward} to develop quasi-Newton methods for the nonsmooth and nonconvex problem \eqref{eq:P3}.
	The gradient interpretation of splitting schemes originated in \cite{rockafellar1976monotone} with the proximal point algorithm and has recently been extended to several other schemes \cite{patrinos2013proximal,patrinos2014douglas,stella2018newton,giselsson2018envelope}.
	In this work we undertake a converse direction: first we design a smooth surrogate of the nonsmooth DC function in \eqref{eq:P}, and then derive a novel splitting algorithm from its gradient steps.
	Classical methods stemming from smooth minimization such as L-BFGS can conveniently be implemented, resulting in a method inherently robust against ill conditioning.
	
	\begin{algorithm}[t]
		\caption{Two-prox algorithm for the DC problem \eqref{eq:P}}%
		\label{alg:P}%
	Select \(\gamma>0\) and \(0<\lambda<2\), and starting from \(s\in\R^p\), repeat
	\begin{equation}\label{eq:alg}
		\begin{cases}[l]
			\left.
				\begin{array}{c @{{}={}} l}
					\fillwidthof[c]{s^+}{u} & \prox_{\gamma\h}(s)
				\\
					\fillwidthof[c]{s^+}{v} & \prox_{\gamma\g}(s)
				\end{array}
			~~\right]
			~~
			\text{\small(in parallel)}
		\\[8pt]
			s^+ = s+\lambda(v-u)
		\end{cases}
	\end{equation}
	{\small
		{\bf Note:}
		\(s^+=s-\lambda\gamma\nabla\env(s)\), where \(\env=\g^\gamma-\h^\gamma\)%
	}%
	\end{algorithm}
	
	\begin{algorithm}[t]
		\caption{Three-prox algorithm for the DC problem \eqref{eq:P3}}
		\label{alg:P3}%
	Select \(0<\gamma<1<\delta\), \(0<\lambda<2(1-\gamma)\), and \(0<\mu<2(1-\delta^{-1})\), and starting from \(s,t\in\R^p\), repeat
	\begin{equation}\label{eq:alg3}
		\begin{cases}[l@{~~}l]
			\left.
				\begin{array}{c @{{}={}} l}
					\fillwidthof[c]{s^+}{u} & \prox_{\frac{\gamma\delta}{\delta-\gamma}\h}\bigl(
						\frac{\delta s-\gamma t}{\delta-\gamma}
					\bigr)
				\\
					\fillwidthof[c]{s^+}{v} & \prox_{\gamma\g}(s)
				\\
					\fillwidthof[c]{s^+}{z} & \prox_{\delta\f}(t)
				\end{array}
			~~\right]
		&
			\text{\small(in parallel)}
		\\[15pt]
			\left.
				\begin{array}{c @{{}={}} l}
					s^+ & s+\lambda(v-u)
				\\
					t^+ & \mathrlap{t+\mu(u-z)}
						\hphantom{\prox_{\frac{\gamma\delta}{\delta-\gamma}\h}\bigl(
						\frac{\delta s-\gamma t}{\delta-\gamma}
					\bigr)}
				\end{array}
			~~\right]
		&
			\text{\small(in parallel)}
		\end{cases}
	\end{equation}
	{\small
		{\bf Note:}
		\begin{tabular}[t]{@{}l@{}}
			\(
				\binom{s^+}{t^+}
			{}={}
				\binom st
				{}-{}
				\binom{\gamma\lambda\I~~\phantom{\delta\mu\I}}{\phantom{\gamma\lambda\I}~~\delta\mu\I}
				\nabla\Psi(s,t)
			\),~~
			where
		\\
			\(\displaystyle
				\Psi(s,t)
			{}={}
				\g^\gamma(s)
				{}-{}
				\f^\delta(t)
				{}-{}
				\h^{\frac{\gamma\delta}{\delta-\gamma}}\bigl(\tfrac{\delta s-\gamma t}{\delta-\gamma}\bigr)
				{}+{}
				\tfrac{1}{2(\delta-\gamma)}\|s-t\|^2
			\)
		\end{tabular}
	}%
	\end{algorithm}

		\subsection{Contributions}
	\paragraph{Fully parallelizable splitting schemes}
		In this paper we propose the novel (sub)gradient-free proximal \cref{alg:P} for the DC problem \eqref{eq:P}, and its fully parallelizable variant when applied to \eqref{eq:P3} synopsized in \cref{alg:P3} (see \cref{sec:Preliminaries} for the notation therein adopted).
		Our approach can be considered complementary to that in \cite{liu2017further}.
		First, we propose a novel smooth DC envelope function (DCE) that shares minimizers and stationary points with the original nonsmooth DC function \(\cost\) in \eqref{eq:P}, similarly to the FBE in \cite{liu2017further}.
		Then, we show that a classical gradient descent on the DCE results in a novel (sub)gradient-free proximal algorithm that is particularly amenable to parallel implementations.
		In fact, even when specialized to problem \eqref{eq:P3} it involves operations on the three functions that can be done in parallel, differently from FBS-based approaches that prescribe serial (sub)gradient and proximal evaluations.
		Due to the complications of computing proximal steps in arbitrary metrics, this flexibility comes at the price of not being able to efficiently handle the composition of \(\f\) in \eqref{eq:P3} with arbitrary linear operators, which is instead possible with FBS-based approaches such as \cite{liu2017further,an2017convergence,banert2019general}.
	
	\paragraph{Novel smooth DC reformulation}
		Thanks to the smooth gradient descent interpretation \emph{it is possible to design classical linesearch strategies} to include directions stemming for instance from quasi-Newton methods, \emph{without complicating the first-order algorithmic oracle}.
		In fact, differently from similar FBE-based quasi-Newton techniques in \cite{liu2017further,patrinos2013proximal,stella2017forward}, no second-order derivatives are needed here and we actually allow for fully nonsmooth formulations.
		Moreover, being the difference of convex and Lipschitz-differentiable functions, the proposed envelope reformulation allows for the extension of the boosted DCA \cite{artacho2018accelerating} to arbitrary DC problems.
	
	\paragraph{A convexity-preserving nonlinear scaling of the FBE}
		When function \(\h\) in \eqref{eq:P} is smooth, we show that the DCE coincides with the FBE \cite{patrinos2013proximal,stella2017forward,themelis2018forward} after a nonlinear scaling.
		This change of variable overcomes some limitations of the FBE, such as preserving convexity when problem \eqref{eq:P} is convex and being (Lipschitz) differentiable without additional requirements on function \(\h\).

		\subsection{Paper organization}
	The paper is organized as follows.
	\Cref{sec:Preliminaries} lists the adopted notational conventions and some known facts needed in the sequel.
	\Cref{sec:env} introduces the DCE, a new envelope function for problem \eqref{eq:P}, and provides some of its basic properties and its connections with the FBE.
	\Cref{sec:Algorithm} shows that a classical gradient method on the DCE results in \cref{alg:P}, and establishes convergence results as a simple byproduct.
	\Cref{alg:P3} is shown to be a \emph{scaled} version of the parent \cref{alg:P}; for the sake of simplicity of presentation, some technicalities needed for this derivation are confined to this section.
	\Cref{sec:Simulations} shows the effect of L-BFGS acceleration on the proposed method on a sparse principal component analysis problem.
	\Cref{sec:Conclusions} concludes the paper.

	\section{Notation and known facts}\label{sec:Preliminaries}
	The set of symmetric matrices in \(\R^p\) is denoted as \(\symm(\R^p)\); the subsets of those which are positive definite is denoted as \(\symm_{++}(\R^p)\).
	Any \(M\in\symm_{++}(\R^p)\) induces the scalar product \((x,y)\mapsto \trans xMy\) on \(\R^p\), with corresponding norm \(\|x\|_M=\sqrt{\trans xMx}\).
	When \(M=\I\), the identity matrix of suitable size, we will simply write \(\|x\|\).
	\(\id\) is the identity function on a suitable space.
	The subdifferential of a proper, lsc, convex function \(\func{\genf}{\R^p}{\Rinf\coloneqq\R\cup\set\infty}\) is
	\[
		\partial\genf(x)
	{}={}
		\set{v\in\R^p}[
			\genf(z)
		{}\geq{}
			\genf(x)
			{}+{}
			\innprod{v}{z-x},~
			\forall z
		].
	\]
	The \DEF{effective domain} of \(\genf\) is \(\dom\genf=\set{x\in\R^p}[\genf(x)<\infty]\), while
	\(
		\conj{\genf}(y)
	\coloneqq{}
		\sup_{x\in\R^p}\set{\innprod xy - \genf(x)}
	\)
	denotes the \DEF{Fenchel conjugate} of \(\f\), which is also proper, closed and convex.
	Properties of conjugate functions are well described for example in \cite{rockafellar1970convex,hiriarturruty2012fundamentals,bauschke2017convex}.
	Among these we recall that
	\begin{equation}\label{eq:ConjSubgr}
		y\in\partial{\genf}(x)
	{}\Leftrightarrow{}
		\innprod xy
		{}={}
		\genf(x)+\conj{\genf}(y)
	{}\Leftrightarrow{}
		x\in\partial\conj{\genf}(y).
	\end{equation}
	The \DEF{proximal mapping} of \(\genf\) with stepsize \(\gamma>0\) is
	\begin{align}\label{eq:prox}
		\prox_{\gamma\genf}(x)
	{}\coloneqq{} &
		\argmin_{w\in\R^p}{
			\set{\genf(w)+\tfrac{1}{2\gamma}\|w-x\|^2}
		},
	\shortintertext{%
		while the value function of the above optimization problem defines the \DEF{Moreau envelope}
	}
		\genf^\gamma(x)
	{}\coloneqq{} &
		\inf_{w\in\R^p}\set{\genf(w)+\tfrac{1}{2\gamma}\|w-x\|^2}.
	\end{align}
	Properties of the Moreau envelope and the proximal mapping are well documented in the literature \cite{bauschke2017convex,
	combettes2005signal,combettes2011proximal}, some of which are summarized next.
	
	\begin{fact}[Proximal properties of convex functions]\label{thm:proxg}%
		Let \(\genf\) be proper, convex, and lsc.
		Then, for all \(\gamma>0\) and \(s,s'\in\R^p\)
		\begin{enumerate}
		\item\label{thm:proxgEquiv}
			\(\prox_{\gamma\genf}(s)\) is the unique point \(x\) such that
			\(
				s\in x+\gamma\partial\genf(x)
			\).
		\item\label{thm:proxgInnprod}
			\(
					\|x-x'\|^2
			{}\leq{}
				\innprod{x-x'}{s-s'}
			{}\leq{}
				\|s-s'\|^2
			\),
			where \(x=\prox_{\gamma\genf}(s)\) and \(x'=\prox_{\gamma\genf}(s')\).
		\item\label{thm:proxgBounds}
			for \(x=\prox_{\gamma\genf}(s)\) and \(w\in\R^p\) it holds that
			\(
				\genf^\gamma(s)
			{}\leq{}
				\genf(w)
				{}+{}
				\tfrac{1}{2\gamma}\|w-s\|^2
				{}-{}
				\tfrac{1}{2\gamma}\|x-s\|^2
			\).
		\item\label{thm:MoreauC1}
			the Moreau envelope \(\genf^\gamma\) is convex and has \(\frac1\gamma\)-Lipschitz-continuous gradient
			\(
				\nabla\genf^\gamma
			{}={}
				\frac1\gamma\bigl(\id-\prox_{\gamma\genf}\bigr)
			\).
		\end{enumerate}
	\end{fact}

	\section{The DC envelope}\label{sec:env}
	In this section we introduce a smooth DC reformulation of \eqref{eq:P} that enables us to cast the nonsmooth and possibly extended-real valued DC problem into the unconstrained minimization of the DCE, a function with Lipschitz-continuous gradient.
	A classical gradient descent algorithm on this reformulation will then be shown in \Cref{sec:Algorithm} to lead to the proposed \cref{alg:P,alg:P3}.
	In this sense, the DCE serves a similar role as the Moreau envelope for the proximal point algorithm \cite{rockafellar1976monotone}, and the FBE and Douglas-Rachford envelope respectively for FBS and the Douglas-Rachford splitting (DRS) \cite{stella2017forward,patrinos2014douglas}.
	
	We begin by formalizing the DC setting of problem \eqref{eq:P} dealt in the paper with the following list of requirements.
	
	\begin{ass}\label{ass:P}%
		The following hold in problem \eqref{eq:P}:
		\begin{enumeratass}
		\item
			\(\func{\g,\h}{\R^p}{\Rinf}\) are proper, convex, and lsc;
		\item\label{ass:phi}
			\(\cost\) is lower bounded (with the convention \(\infty-\infty=\infty\)).
		\end{enumeratass}
	\end{ass}
	
	\begin{defin}[DC envelope]\label{defin:env}
		Suppose that \cref{ass:P} holds.
		Relative to problem \eqref{eq:P}, the DC envelope (DCE) with stepsize \(\gamma>0\) is
		\[
			\env(s)
		{}\coloneqq{}
			\g^\gamma(s)-\h^\gamma(s).
		\]
	\end{defin}
	
	Before showing that the DCE \(\env\) satisfies the anticipated smoothness properties and is tightly connected with solutions of problem \eqref{eq:P}, we provide a simple characterization of stationary points in terms of the proximal mappings of the functions involved in the DC formulation.
	This will then be used to connect points that are stationary in the sense of \eqref{eq:stationary} for \eqref{eq:P} with points that are stationary in the classical sense for \(\env\).
	
	\begin{lem}[Optimality conditions]\label{thm:optimality}%
		Suppose that \cref{ass:P} holds.
		Then, any of the following is equivalent to stationarity at \(u\) in the sense of \eqref{eq:stationary}:
		\begin{enumerateq}
		\item\label{thm:prox=}
			there exist \(\gamma>0\) and \(s\in\R^p\) such that \(u=\prox_{\gamma\g}(s)=\prox_{\gamma\h}(s)\);
		\item\label{thm:prox=all}
			for all \(\gamma>0\) there exists \(s\in\R^p\) such that \(u=\prox_{\gamma\g}(s)=\prox_{\gamma\h}(s)\).
		\end{enumerateq}
		\begin{proof}
			If \(u\) is stationary, then for every \(\gamma>0\) and \(\xi\in\partial\g(u)\cap\partial\h(u)\neq\emptyset\) it follows from \cref{thm:proxgEquiv} that
			\(
				u=\prox_{\gamma\g}(s)=\prox_{\gamma\g}(s)
			\)
			for \(s=u+\gamma\xi\), proving \ref{thm:prox=all} and thus \ref{thm:prox=}.
			Conversely, if \ref{thm:prox=} holds then \cref{thm:proxgEquiv} again implies
			\(
				\frac{s-u}{\gamma}\in\partial g(u)
			\)
			and
			\(
				\frac{s-u}{\gamma}\in\partial h(u)
			\),
			proving that \(u\) is stationary.
		\end{proof}
	\end{lem}
	
	\begin{lem}[Basic properties of the DCE]\label{thm:env}
		Let \cref{ass:P} hold, and for notational conciseness given \(s\in\R^p\) let \(u\coloneqq\prox_{\gamma\h}(s)\) and \(v\coloneqq\prox_{\gamma\g}(s)\).
		The following hold:
		\begin{enumerate}
		\item\label{thm:smooth}
			\(\env\) is \(\tfrac1\gamma\)-smooth with
			\(
				\nabla\env=\tfrac1\gamma\bigl(\prox_{\gamma\h}-\prox_{\gamma\g}\bigr)
			\);
		\item\label{thm:stationary}
			\(\nabla\env(s)=0\) iff \(u\) is stationary (cf. \eqref{eq:stationary});
		\item\label{thm:sandwich}
			\(
				\cost(v)
				{}+{}
				\tfrac{1}{2\gamma}\|v-u\|^2
			{}\leq{}
				\env(s)
			{}\leq{}
				\cost(u)
				{}-{}
				\tfrac{1}{2\gamma}\|v-u\|^2
			\);
		\item\label{thm:min}
			\(
				\argmin\cost
			{}={}
				\prox_{\gamma\h}(S_\star)
			{}={}
				\prox_{\gamma\g}(S_\star)
			\)
			and
			\(\inf\cost=\inf\env\)
			for \(S_\star=\argmin\env\).
		\end{enumerate}
		\begin{proof}~
			\begin{proofitemize}
			\item\ref{thm:smooth}~
				The expression of the gradient follows from \cref{thm:MoreauC1}.
				The bounds in \cref{thm:proxgInnprod} imply that
				\begin{equation}\label{eq:smoothbounds}
					\left|\innprod{\nabla\env(s)-\nabla\env(s')}{s-s'}\right|
				{}\leq{}
					\tfrac1\gamma
					\|s-s'\|^2,
				\end{equation}
				proving that \(\nabla\env\) is \(\gamma^{-1}\)-Lipschitz continuous.
			\item\ref{thm:stationary}
				Follows from assertion \ref{thm:smooth} and \cref{thm:optimality}.
			\item\ref{thm:sandwich}
				Follows by applying the proximal inequalities of \cref{thm:proxgBounds} with \(w=u\) and \(w=v\).
			\item\ref{thm:min}
				Follows from assertion \ref{thm:sandwich}, \cref{thm:optimality}, and the fact that global minimizers for \(\cost\) are stationary.
			\qedhere
			\end{proofitemize}
		\end{proof}
	\end{lem}

		\subsection{Connections with the forward-backward envelope}
	As it will be detailed in \Cref{sec:hypo}, considering difference of hypo\-convex functions in problem \eqref{eq:P} leads to virtually no generalization.
	A more interesting scenario occurs when both \(\h\) and \(-\h\) are hypoconvex functions, which amounts to \(\h\) being \(L_{\h}\)-smooth (differentiable with \(L_{\h}\)-Lipschitz gradient).
	In order to elaborate this property we first need to specialize \cref{thm:proxf} to smooth functions.
	
	\begin{lem}[Proximal properties of smooth functions]\label{thm:proxf}%
		Suppose that \(\func{\genf}{\R^p}{\R}\) is \(L_{\genf}\)-smooth.
		Then, there exist
		\(
			\sigma_{\genf},\sigma_{-\genf}
		{}\in{}
			[-L_{\genf},L_{\genf}]
		\)
		with
		\(
			L_{\genf}
		{}={}
			\max\set{
				|\sigma_{\genf}|,|\sigma_{-\genf}|
			}
		\)
		such that \(\genf-\tfrac{\sigma_{\genf}}{2}\|{}\cdot{}\|^2\) and \(-\genf-\tfrac{\sigma_{-\genf}}{2}\|{}\cdot{}\|^2\) are convex functions.
		Then, for all \(\gamma<\nicefrac{1}{[\sigma_{-\genf}]_-}\) (with the convention \(\nicefrac10=\infty\)) and \(s,s'\in\R^p\)
		\begin{enumerate}
		\item\label{thm:proxfEquiv}
			\(\prox_{-\gamma\genf}(s)\) is the unique \(u\) such that
			\(
				s=u-\gamma\nabla\genf(u)
			\);%
		\item\label{thm:proxfInnprod}
			\(
				\tfrac{1}{1-\gamma\sigma_{\genf}}
				\|s-s'\|^2
			{}\leq{}
				\innprod{u-u'}{s-s'}
			{}\leq{}
				\tfrac{1}{1+\gamma\sigma_{-\genf}}
				\|s-s'\|^2
			\),
			where \(u=\prox_{-\gamma\genf}(s)\) and \(u'=\prox_{-\gamma\genf}(s')\);
		\item
			\((-\genf)^\gamma\) is differentiable with \(\nabla(-\genf)^\gamma=\frac{\id-\prox_{-\gamma\genf}}{\gamma}\).
		\end{enumerate}
		\begin{proof}
			The claim on the existence of \(\sigma_{\pm\genf}\) comes from the fact that \(\genf\) is \(L_{\genf}\)-smooth iff \(\tfrac{L_{\genf}}{2}\|{}\cdot{}\|^2\pm\genf\) are convex functions, and that \(\genf\) is \(L_{\genf}\)-smooth iff so is \(-\genf\).
			All other claims then follow from \cref{thm:proxg} applied to the convex function \(\tilde\genf=-\genf-\tfrac{\sigma_{-\genf}}{2}\|{}\cdot{}\|^2\), in light of the identity
			\(
				\prox_{\gamma\tilde\genf}
			{}={}
				\prox_{-\frac{\gamma}{1-\gamma\sigma_{-\genf}}\genf}
				{}\circ{}
				\tfrac{\id}{1-\gamma\sigma_{-\genf}}
			\)
			\cite[Prop. 24.8(i)]{bauschke2017convex}.
		\end{proof}
	\end{lem}
	
	In the remainder of this subsection, suppose that \(\h\) is smooth.
	Denoting \(\f\coloneqq-\h\), problem \eqref{eq:P} reduces to
	\begin{equation}\label{eq:FB:P}
		\minimize_{u\in\R^n}\f(u)+\g(u)=\g(u)-(-\f)(u)
	\end{equation}
	with \(\g\) convex and \(\f\) smooth.
	A textbook algorithm for addressing such composite minimization problems is FBS, which interleaves proximal and gradient operations as
	\begin{equation}
		u^+=\FB u.
	\end{equation}
	By observing that \(s=\Fw u\) iff \(u=\prox_{-\gamma\f}(s)\) for \(\gamma<\nicefrac{1}{L_f}\), one obtains the following curious connection among \(\env\) and the forward-backward envelope \cite[Eq. (2.3)]{stella2017forward}
	\begin{equation}\label{eq:FBE}
		\FBE(u)
	{}={}
		\f(u)
		{}-{}
		\tfrac\gamma2\|\nabla\f(u)\|^2
		{}+{}
		\g^\gamma(\Fw u).
	\end{equation}
	
	\begin{lem}
		\renewcommand{\envh}{-\f}%
		In problem \eqref{eq:FB:P}, suppose that \(\f\) is \(L_{\f}\)-smooth and \(\g\) is proper, convex, and lsc.
		Then, for every \(\gamma<\nicefrac{1}{L_{\f}}\)
		\[
			\FBE
		{}={}
			\env\circ(\Fw{})
		~\text{and}~
			\env
		{}={}
			\FBE\circ\prox_{-\gamma\f}.
		\]
		Moreover, \(\env\) is \(\frac{1-\gamma L_{\f}}{\gamma}\)-smooth, and if \(\f\) is additionally convex then so is \(\env\).
		\begin{proof}
			Let \(u\in\R^p\) and \(\gamma\in(0,\nicefrac{1}{L_{\f}})\) be fixed, and for notational conciseness let \(u=\prox_{-\gamma\f}(s)\).
			Then, \(s=\Fw u\) and \((-\f)^\gamma(s)=-\f(u)+\tfrac{1}{2\gamma}\|u-s\|^2\), hence
			\begin{align*}
				\env(s)
			{}={} &
				\g^\gamma(\Fw u)
				{}+{}
				\f(u)-\tfrac{1}{2\gamma}\|u-s\|^2
			\\
			{}={} &
				\f(u)
				{}-{}
				\tfrac\gamma2\|\nabla\f(u)\|^2
				{}+{}
				\g^\gamma(\Fw u),
			\end{align*}
			which is exactly \(\FBE(u)\), cf. \eqref{eq:FBE}.
			By using \cref{thm:proxfInnprod} for \(\h=-\f\), the bounds in \eqref{eq:smoothbounds} become
			\[
				\tfrac{\sigma_{\f}\|s-s'\|^2}{1-\gamma\sigma_{\f}}
			{}\leq{}
				\innprod{\nabla\env(s)-\nabla\env(s')}{s-s'}
			{}\leq{}
				\tfrac{\gamma^{-1}\|s-s'\|^2}{1+\gamma\sigma_{-\f}}.
			\]
			Since \(|\sigma_{\f}|,|\sigma_{-\f}|\leq L_{\f}\), the claimed smoothness follows.
			Finally, if \(\f\) is convex then \(\sigma_{\f}\) is nonnegative and thus so is the lower bound above, proving convexity of \(\env\).
		\end{proof}
	\end{lem}


	\section{The algorithm}\label{sec:Algorithm}
	Having assessed the \(\frac1\gamma\)-smoothness of \(\env\) and its connection with problem \eqref{eq:P} in \cref{thm:env}, the minimization of the nonsmooth DC function \(\cost=\g-\h\) can be carried out with a gradient descent with constant stepsize \(\tau<2\gamma\) on \(\env\).
	As shown in the next result, this is precisely \cref{alg:P}.
	
	\begin{thm}\label{thm:alg}%
		Suppose that \cref{ass:P} holds, and starting from \(s^0\in\R^n\) consider the iterates \(\seq{s^k,u^k,v^k}\) generated by \cref{alg:P} with \(\gamma>0\) and \(\lambda\in(0,2)\).
		Then, for every \(k\in\N\) it holds that
		\(
			s^{k+1}
		{}={}
			s^k
			{}-{}
			\gamma\lambda\nabla\env(s^k)
		\)
		and
		\begin{equation}\label{eq:SD}
			\env(s^{k+1})
		{}\leq{}
			\env(s^k)
			{}-{}
			\tfrac{\lambda(2-\lambda)}{2\gamma}
			\|u^k-v^k\|^2.
		\end{equation}
		In particular:
		\begin{enumerate}
		\item\label{thm:alg:res}
			the fixed-point residual vanishes with \(\min_{i\leq k}\|u^i-v^i\|=o(\nicefrac{1}{\sqrt k})\);
		\item\label{thm:alg:omega}%
			\(\seq{u^k}\) and \(\seq{v^k}\) have the same set of cluster points, be it \(\Omega\); when \(\seq{s^k}\) is bounded, every \(u_\star\in\Omega\) is stationary for \(\cost\) (in the sense of \eqref{eq:stationary}) and \(\cost\) is constant on \(\Omega\), the value being the (finite) limit of the sequences \(\seq{\env(s^k)}\) and \(\seq{\cost(v^k)}\);%
		\item\label{thm:alg:bounded}
			if \(\cost\) is coercive, then \(\seq{s^k,u^k,v^k}\) is bounded.
		\end{enumerate}
		\begin{proof}
			That \(s^{k+1}=s^k-\lambda\gamma\nabla\env(s^k)\) follows from \cref{thm:smooth}.
			The proof is now standard, see \eg \cite{bertsekas2016nonlinear}:
			\(\frac1\gamma\)-smoothness implies the upper bound
			\begin{align*}
				\env(s^{k+1})
			{}\leq{} &
				\env(s^k)
				{}+{}
				\innprod{\nabla\env(s^k)}{s^{k+1}-s^k}
			\\ &
				{}+{}
				\tfrac{1}{2\gamma}\|s^{k+1}-s^k\|^2
			\\
			{}={} &
				\env(s^k)
				{}-{}
				\tfrac{\lambda(2-\lambda)}{2\gamma}
				\|u^k-v^k\|^2,
			\end{align*}
			which is \eqref{eq:SD}.
			We now show the numbered claims.
			\begin{proofitemize}
			\item\ref{thm:alg:res}
				By telescoping \eqref{eq:SD} and using the fact that \(\inf\env=\inf\cost>-\infty\) owing to \cref{thm:min,ass:phi}, we obtain that the sequence of squared residuals \(\seq{\|u^k-v^k\|^2}\) has finte sum, hence the claim.
			\item\ref{thm:alg:omega}
				That the sequences have same cluster points follows from assertion \ref{thm:alg:res}.
				Moreover, \eqref{eq:SD} and the lower boundedness of \(\env\) imply that the sequence \(\seq{\env(s^k)}\) monotonically decreases to a finite value, be it \(\cost_\star\).
				Continuity of \(\env\) then implies that \(\env(s_\star)=\varphi_\star\) for every limit point \(s_\star\) of \(\seq{s^k}\).
				If \(\seq{s^k}\) is bounded, then so are \(\seq{u^k}\) and \(\seq{v^k}\) owing to Lipschitz continuity of the proximal mappings.
				Moreover, for every \(k\) one has \(s^k=u^k+\gamma\xi^k=v^k+\gamma\eta^k\) for some \(\xi^k\in\partial\h(u^k)\) and \(\eta^k\in\partial\g(v^k)\).
				Necessarily, the sequences of subgradients are bounded, and for any limit point \(u_\star\) of \(\seq{u^k}\), up to possibly extracting, we have that \(u_\star=\prox_{\gamma\h}(s_\star)=\prox_{\gamma\g}(s_\star)\) for some cluster point \(s_\star\) of \(\seq{s^k}\).
				By invoking \cref{thm:optimality} we conclude that \(\cost(u_\star)=\varphi_\star\).
			\item\ref{thm:alg:bounded}
				Boundedness of \(\seq{s^k}\) follows from the fact that \(\env(s^k)\leq\env(s^0)\) for all \(k\), owing to \eqref{eq:SD}.
				In turn, boundedness of \(\seq{u^k}\) and \(\seq{v^k}\) follows from Lipschitz continuity of the proximal mappings.
			\qedhere
			\end{proofitemize}
		\end{proof}
	\end{thm}
	
	The remainder of the section is devoted to deriving \cref{alg:P3} as a special instance of \cref{alg:P} applied to the problem reformulation \eqref{eq:P3DC}.
	In order to formalize this derivation, we first need to address a minor technicality arising because of the nonconvexity of function \(\H\) therein, which prevents a direct application of \cref{alg:P} to the function decomposition \(\G-\H\).
	Fortunately however, by simply adding a quadratic term to both \(\G\) and \(\H\) the desired DC formulation is obtained without actually changing the cost function \(\Cost\) in problem \eqref{eq:P3DC}.
	This simple issue is addressed next.

		\subsection{Strongly and hypoconvex functions}\label{sec:hypo}
	Clearly, adding a same quantity to both functions \(\g\) and \(\h\) leaves problem \eqref{eq:P} unchanged.
	In particular, the convexity setting of \cref{ass:P} can also be achieved when \(\g\) and \(\h\) are \emph{hypoconvex}, in the sense that they are convex up to adding a suitably large quadratic function.
	Recall that for \(\tilde\genf=\genf+\tfrac\mu2\|{}\cdot{}\|^2\) it holds that
	\(
		\prox_{\tilde\gamma\tilde\genf}(s)
	{}={}
		\prox_{\gamma\genf}(\tfrac{s}{1+\tilde\gamma\mu})
	\)
	for \(\gamma=\frac{\tilde\gamma}{1+\tilde\gamma\mu}\) \cite[Prop. 24.8(i)]{bauschke2017convex}.
	Therefore, as long as there exists \(\mu\in\R\) such that both \(\g+\tfrac\mu2\|{}\cdot{}\|^2\) and \(\h+\tfrac\mu2\|{}\cdot{}\|^2\) are convex functions, one can apply iterations \eqref{eq:alg} to the minimization of
	\(
		\g+\tfrac\mu2\|{}\cdot{}\|^2
	{}-{}
		\left(
			\h+\tfrac\mu2\|{}\cdot{}\|^2
		\right)
	\)
	to obtain
	\[
		\begin{cases}[l @{{}={}} l]
			u^k & \prox_{\tilde\gamma\h}(\tilde s^k)
		\\
			v^k & \prox_{\tilde\gamma\g}(\tilde s^k)
		\\[2pt]
			\tilde s^{k+1} & \tilde s^k+\tilde\lambda(v^k-u^k),
		\end{cases}
	\]
	where \(\tilde\gamma\coloneqq\frac{\gamma}{1+\gamma\mu}\), \(\tilde s^k\coloneqq\frac{1}{1+\gamma\mu}s^k\), and \(\tilde\lambda\coloneqq\frac{1}{1+\gamma\mu}\lambda\).
	By observing that
	\(
		\frac{\gamma}{1+\gamma\mu}
	\)
	ranges in \((0,\nicefrac1\mu)\) for \(\gamma\in(0,\infty)\) (with the convention \(\nicefrac10=\infty\)), and that
	\(
		\tilde\lambda
	{}={}
		\lambda(1-\tilde\gamma\mu)
	\),
	we obtain the following.
		
	\begin{rem}[\emph{Strongly} convex and \emph{hypo}convex functions]\label{thm:hypo}%
		If \(\mu\in\R\) is such that both \(\g+\tfrac\mu2\|{}\cdot{}\|^2\) and \(\h+\tfrac\mu2\|{}\cdot{}\|^2\) are convex functions, then all the numbered claims of \cref{thm:alg} still hold provided that \(0<\lambda<2(1-\gamma\mu)\).
	\end{rem}
	
	As a final step towards the analysis of \cref{alg:P3}, in the next subsection we motivate the presence of the two additional parameters \(\delta\) and \(\mu\) missing in \cref{alg:P}.

		\subsection{Matrix stepsize and relaxation}
	A substantial degree of flexibility can be introduced by replacing the quadratic term \(\tfrac{1}{2\gamma}\|w-{}\cdot{}\|^2\) appearing in the definition \eqref{eq:prox} of the proximal mapping with the squared norm \(\tfrac12\|w-{}\cdot{}\|_{\Gamma^{-1}}^2\) induced by a matrix \(\Gamma\in\symm_{++}(\R^p)\).
	The scalar stepsize \(\gamma\) is achieved by considering \(\Gamma=\gamma\I\); in general, we may thus think of \(\Gamma\) as a matrix stepsize.
	Denoting
	\begin{align}\label{eq:PROX}
		\prox_{\genf}^\Gamma(x)
	{}={} &
		\argmin_w\set{\genf(w)+\tfrac12\|w-x\|_{\Gamma^{-1}}^2}
	\shortintertext{and}
		\genf^\Gamma(x)
	{}={} &
		\min_w\set{\genf(w)+\tfrac12\|w-x\|_{\Gamma^{-1}}^2}
	\end{align}
	the corresponding Moreau envelope, as shown in \cite[Thm. 4.1.4]{hiriarturruty1993convex} we have that \(\nabla\genf^\Gamma=\Gamma^{-1}(\id-\prox_{\genf}^\Gamma)\) satisfies
	\[
		0
	{}\leq{}
		\innprod{\nabla\genf^\Gamma(s)-\nabla\genf^\Gamma(s')}{s-s'}
	{}\leq{}
		\|s-s'\|_{\Gamma^{-1}}^2.
	\]
	
	\begin{rem}[Matrix stepsizes and relaxations]\label{thm:matrix}%
		\renewcommand{\envgamma}{\Gamma}%
		Under \cref{ass:P}, given a diagonal stepsize \(\Gamma\in\symm_{++}(\R^p)\) and a diagonal relaxation \(\Lambda\in\symm_{++}(\R^p)\) the iterations
		\begin{equation}\label{eq:ALG}
			\begin{cases}[l @{{}={}} l]
				u^k & \prox_{\h}^\Gamma(s^k)
			\\[2pt]
				v^k & \prox_{\g}^\Gamma(s^k)
			\\[2pt]
				s^{k+1} & s^k+\Lambda(v^k-u^k)
			\end{cases}
		\end{equation}
		produce a sequence such that
		\[
			\env(s^{k+1})
		{}\leq{}
			\env(s^k)
			{}-{}
			\tfrac12
			\|u^k-v^k\|_{(2\I-\Lambda)\Gamma^{-1}\Lambda}^2.
		\]
		In particular, all the numbered claims of \cref{thm:alg} still hold when
		\(
			0\prec\Lambda\prec 2\I
		\).\footnote{%
			Although similar claims can be made for more general positive definite matrices, the diagonal requirement guarantees the symmetry of \((2\I-\Lambda)\Gamma^{-1}\Lambda\) and thus its positive definiteness for \(\Lambda\) as prescribed above.%
		}%
	\end{rem}
	
	Notice that the optimality condition for minimization problem \eqref{eq:PROX} reads
	\(
		0
	{}\in{}
		\partial\genf(w)
		{}+{}
		\Gamma^{-1}(w-x)
	\).
	Equivalently,
	\begin{equation}\label{eq:PROXequiv}
		w=\prox_{\genf}^\Gamma(x)
	\quad\Leftrightarrow\quad
		x\in w+\Gamma\partial\genf(w).
	\end{equation}
	By using this fact, if a symmetric matrix \(M\) is such that the function \(\tilde\genf=\genf+\tfrac12\innprod{{}\cdot{}}{M{}\cdot{}}\) is convex, one can express its proximal map in terms of that of \(\genf\) in a similar fashion as the scalar case considered in \cref{sec:hypo}, namely,
	\[
		\prox_{\tilde\genf}^{\tilde\Gamma}
	{}={}
		\prox_{\genf}^{\Gamma}\circ(\I-\Gamma M)
	\]
	with \(\Gamma=(\tilde\Gamma^{-1}+M)^{-1}\).\footnote{%
		These expressions in terms of the new stepsize \(\Gamma\) use the matrix identities
		\(
			(\I+\tilde\Gamma M)^{-1}\tilde\Gamma
		{}={}
			(\tilde\Gamma^{-1}+M)^{-1}
		\)
		and
		\(
			(\I+\tilde\Gamma M)^{-1}
		{}={}
			\I-\Gamma M
		\)
		for
		\(
			\Gamma=(\I+\tilde\Gamma M)^{-1}\tilde\Gamma
		\).
	}
	It is thus possible to combine \cref{thm:hypo,thm:matrix} as follows, where again for simplicity we restrict the case to diagonal matrices.
	
	\begin{rem}\label{thm:hypomatrix}%
		\renewcommand{\envgamma}{\Gamma}%
		If a diagonal matrix \(M\) is such that both functions \(\g+\tfrac12\innprod{{}\cdot{}}{M{}\cdot{}}\) and \(\h+\tfrac12\innprod{{}\cdot{}}{M{}\cdot{}}\) are convex, then the sequence produced by \eqref{eq:ALG} satisfies all the numbered claims of \cref{thm:alg} as long as
		\(
			0
		{}\prec{}
			\Lambda
		{}\prec{}
			2(\I-\Gamma M)
		\).
	\end{rem}

		\subsection{A parallel three-prox splitting}\label{sec:3splitting}
	After the generalization documented in \cref{thm:hypomatrix} we are ready to address the formulation \eqref{eq:P3} and express \cref{alg:P3} as a ``scaled'' variant of \cref{alg:P}.
	We begin by rigorously framing the problem setting.
	
	\begin{ass}\label{ass:P3}
		In problem \eqref{eq:P3}
		\begin{enumeratass}
		\item
			\(\func{f,g,h}{\R^n}{\Rinf}\) are proper, lsc, and convex;
		\item
			\(\cost\) is lower bounded.
		\end{enumeratass}
	\end{ass}
	
	\begin{thm}\label{thm:alg3}%
		\renewcommand{\envg}{\G}%
		\renewcommand{\envh}{\H}%
		\renewcommand{\envgamma}{\Gamma}%
		Let \cref{ass:P3} hold, and starting from \((s^0,t^0)\in\R^n\times\R^n\) consider the iterates \(\seq{s^k,t^k,u^k,v^k,z^k}\) generated by \cref{alg:P3} with \(0<\gamma<1<\delta\), \(0<\lambda<2(1-\gamma)\) and \(0<\mu<2(1-\delta^{-1})\).
		Then, denoting
		\begin{align*}
			\Psi(s,t)
		{}={} &
			\env(s,t\nicefrac{}{\delta})
		\\
		\numberthis\label{eq:Psi3}
		{}={} &
			\g^\gamma(s)
			{}-{}
			\f^\delta(t)
			{}-{}
			\h^{\frac{\gamma\delta}{\delta-\gamma}}\bigl(\tfrac{\delta s-\gamma t}{\delta-\gamma}\bigr)
			{}+{}
			\tfrac{1}{2(\delta-\gamma)}\|s-t\|^2,
		\end{align*}
		for every \(k\in\N\) it holds that
		\begin{equation}\label{eq:P3:GD}
		\textstyle
			\binom{s^{k+1}}{t^{k+1}}
		{}={}
			\binom{s^k}{t^k}
			{}-{}
			\binom{\gamma\lambda\I~~\phantom{\delta\mu\I}}{\phantom{\gamma\lambda\I}~~\delta\mu\I}
			\nabla\Psi(s^k,t^k).
		\end{equation}
		Moreover
		\begin{enumerate}
		\item\label{thm:alg3:res}
			the fixed-point residual vanishes with \(\min_{i\leq k}\|\binom{u^i-v^i}{u^i-z^i}\|=o(\nicefrac{1}{\sqrt k})\);
		\item\label{thm:alg3:omega}%
			\(\seq{u^k}\) \(\seq{v^k}\) and \(\seq{z^k}\) have the same set of cluster points, be it \(\Omega\); when \(\seq{s^k}\) is bounded, every \(u_\star\in\Omega\) satisfies the stationarity condition
			\[
				\emptyset
			{}\neq{}
				\partial\g(u_\star)
				{}\cap{}
				\bigl(
					\partial\f(u_\star)
					{}+{}
					\partial\h(u_\star)
				\bigr)
			{}\subseteq{}
				\partial\g(u_\star)
				{}\cap{}
				\partial(\f+\h)(u_\star)
			\]
			and \(\cost\) is constant on \(\Omega\), the value being the (finite) limit of the sequence \(\seq{\cost(u^k)}\);%
		\item\label{thm:alg3:bounded}
			if \(\cost\) is coercive, then \(\seq{s^k,t^k,u^k,v^k,z^k}\) is bounded.
		\end{enumerate}
		\begin{proof}
			Let \(\Cost\), \(\G\) and \(\H\) be as in \eqref{eq:P3DC}, and observe that
			\[
				\Phi(x,y)
			{}\geq{}
				\inf_{y'}\Phi(x,y')
			{}={}
				\varphi(x).
			\]
			In particular, if \(\cost\) is coercive then necessarily so is \(\Cost\).
			Let
			\(
				\Gamma
			{}\coloneqq{}
				\binom{\gamma\I~~\phantom{\delta^{-1}\I}}{\phantom{\gamma\I}~~\delta^{-1}\I}
			\).
			Under \cref{ass:P3}, function \(\G\) is convex and one can easily verify that
			\begin{align*}
				(v_s,v_t)
			{}={}
				\prox_{\G}^\Gamma(s,t)
			{}\Leftrightarrow{} &
				\begin{cases}
					v_s
				{}={}
					\prox_{\gamma\g}(s)
				\\
					v_t
				{}={}
					t-\delta^{-1}\prox_{\delta\f}(\delta t)
				\end{cases}
			\shortintertext{%
				in light of the Moreau identity
				\(
					\prox_{\nicefrac{\conj{\f}}{\delta}}(t)
				{}={}
					t-\delta^{-1}\prox_{\delta\f}(\delta t)
				\),
				see \cite[Thm. 14.3(ii)]{bauschke2017convex}.
				Furthermore, from \eqref{eq:PROXequiv} we have
			}
				(u_s,u_t)
			{}={}
				\prox_{\H}^\Gamma(s,t)
			{}\Leftrightarrow{}&
				\begin{cases}
					s
				{}\in{}
					u_s+\gamma\partial\h(u_s)+\gamma u_t
				\\
					t
				{}={}
					u_t+u_s\nicefrac{}{\delta}
				\end{cases}
			\\
			{}\Leftrightarrow{}&
				\begin{cases}
					\frac{s-\gamma t}{1-\nicefrac\gamma\delta}
				{}\in{}
					u_s+\frac{\gamma}{1-\nicefrac\gamma\delta}\partial\h(u_s)
				\\
					u_t
				{}={}
					t-u_s\nicefrac{}{\delta}
				\end{cases}
			\\
			{}\Leftrightarrow{}&
				\begin{cases}
					u_s
				{}={}
					\prox_{\frac{\gamma\delta}{\delta-\gamma}\h}\bigl(
						\frac{\delta s-\gamma\delta t}{\delta-\gamma}
					\bigr)
				\\
					u_t
				{}={}
					t-u_s\nicefrac{}{\delta}.
				\end{cases}
			\end{align*}
			In particular,
			\[\textstyle
				\binom{s}{\delta t}
				{}+{}
				\binom{\lambda\I~~\phantom{\delta\mu\I}}{\phantom{\lambda\I}~~\delta\mu\I}
				\left(
					\prox_{\G}^\Gamma\binom st
					{}-{}
					\prox_{\H}^\Gamma\binom st
				\right)
			{}={}
				\binom{
					s+\lambda(v_s-u_s)
				}{
					\delta t+\mu(u_s-\prox_{\delta\f}(\delta t)
				}.
			\]
			Apparently, iterations \eqref{eq:alg3} correspond to those in \eqref{eq:ALG} with
			\(
				\Lambda
			{}\coloneqq{}
				\binom{\lambda\I~~\phantom{\mu\I}}{\phantom{\lambda\I}~~\mu\I}
			\)
			after the scaling \(t\gets\nicefrac t\delta\).
			From these computations and using the fact that
			\(
				(\conj{\f})^{\nicefrac1\delta}\circ\nicefrac{\id}{\delta}
			{}={}
				\tfrac{1}{2\delta}\|{}\cdot{}\|^2
				{}-{}
				\f^\delta
			\),
			see \cite[Thm. 14.3(i)]{bauschke2017convex}, the expressions in \eqref{eq:Psi3} and \eqref{eq:P3:GD} are obtained.
			Since function \(\H+\tfrac12\|{}\cdot{}\|^2\) is convex --- that is, the setting of \cref{thm:hypomatrix} is satisfied with \(M=\I\) --- and the condition
			\(
				0
			{}\prec{}
				\Lambda
			{}\prec{}
				2(\I-\Gamma)
			\)
			holds when \(\gamma,\delta,\lambda,\mu\) are as in the statement, it only remains to show that the limit points satisfy the stationarity condition of assertion \ref{thm:alg3:omega}, as the rest of the proof follows from \cref{thm:alg:res,thm:hypomatrix}.
			To this end, since
			\(
				\binom{v^k-u^k}{u^k-z^k}
			{}={}
				\binom{s^{k+1}-s^k}{t^{k+1}-t^k}
			{}\to{}
				0
			\)
			the sequences \(\seq{u^k}\) \(\seq{v^k}\) and \(\seq{z^k}\) have the same cluster points.
			If \(\seq{s^k}\) is bounded, arguing as in the proof of \cref{thm:alg:omega} we have that if \(u^k\to u_\star\) as \(k\in K\) for an infinite set of indices \(K\subseteq\N\), necessarily also \(v^k\to u_\star\) as \(k\in K\), and \((s^k,t^k)\to (s_\star,t_\star)\) as \(k\in K\) for some \(s_\star,t_\star\) such that
			\[
				\prox_{\frac{\gamma\delta}{\delta-\gamma}\h}
				\bigl(
					\tfrac{\delta s_\star-\gamma t_\star}{\delta-\gamma}
				\bigr)
			{}={}
				\prox_{\gamma\g}(s_\star)
			{}={}
				\prox_{\delta\f}(t_\star).
			\]
			We then conclude from \cref{thm:proxgEquiv} that
			\[
				\frac{
					\frac{\delta s_\star-\gamma t_\star}{\delta-\gamma}
					{}-{}
					u_\star
				}{
					\frac{\gamma\delta}{\delta-\gamma}
				}
			{}\in{}
				\partial\h(u_\star),
			~
				\frac{s_\star-u_\star}{\gamma}
			{}\in{}
				\partial\g(u_\star),
			~
				\frac{t_\star-u_\star}{\delta}
			{}\in{}
				\partial\f(u_\star),
			\]
			which gives
			\[
				\tfrac{s_\star-u_\star}{\gamma}
			{}\in{}
				\partial\g(u_\star)
				{}\cap{}
				\left(
					\partial\f(u_\star)
					{}+{}
					\partial\h(u_\star)
				\right),
			\]
			and the claimed stationarity condition follows from the inclusion
			\(
				\partial\f+\partial\h
			{}\subseteq{}
				\partial(\f+\h)
			\),
			see \cite[Thm. 23.8]{rockafellar1970convex}.
		\end{proof}
	\end{thm}

	\section{Simulations}\label{sec:Simulations}
	We study the performance of \cref{alg:P} applied to a sparse principal component analysis (SPCA) problem.
	Following \cite[\S2.1]{journee2010generalized}, an SPCA problem can be formulated as 
	\begin{equation}\label{eq:PCA}
		\minimize -\tfrac12\trans s\Sigma s + \kappa\|s\|_1
	\quad\stt{}
		s\in\cball01
	\end{equation}
	with \(\cball01\coloneqq\set{s}[\|s\|\leq 1]\), \(\Sigma = \trans AA\) the sample covariance matrix, and \(\kappa\) a sparsity inducing parameter.
	This problem can be identified as a DC problem of type \eqref{eq:P} by denoting \(\g(s) = \kappa\|s\|_1 + \indicator_{\cball01}(s)\) and \(\h(s) = \tfrac{1}{2}\trans s\Sigma s\), where \(\indicator_C\) denotes the indicator function of a (nonempty closed convex) set \(C\), namely \(\indicator_C(x)=0\) if \(x\in C\) and \(\infty\) otherwise.
	Then,
	\begin{align*}
		\prox_{\gamma\h}(s)
	{}={} &
		(\I +\gamma\Sigma)^{-1}s,
	~~\text{and}
	\\
		\prox_{\gamma\g}(s)
	{}={} &
		\frac{\sign(s)\odot[|s|-\kappa\gamma{\bf 1}]_+}{\max\set{1,\|[|s|-\kappa\gamma{\bf 1}]_+\|}},
	\end{align*}
	with \(\odot\) the elementwise multiplication, \(|{}\cdot{}|\) the elementwise absolute value, and \({\bf 1}\) the \(\R^n\)-vector of all ones. 
	
\begin{figure}
		\centering
		\includetikz[width=8.0cm]{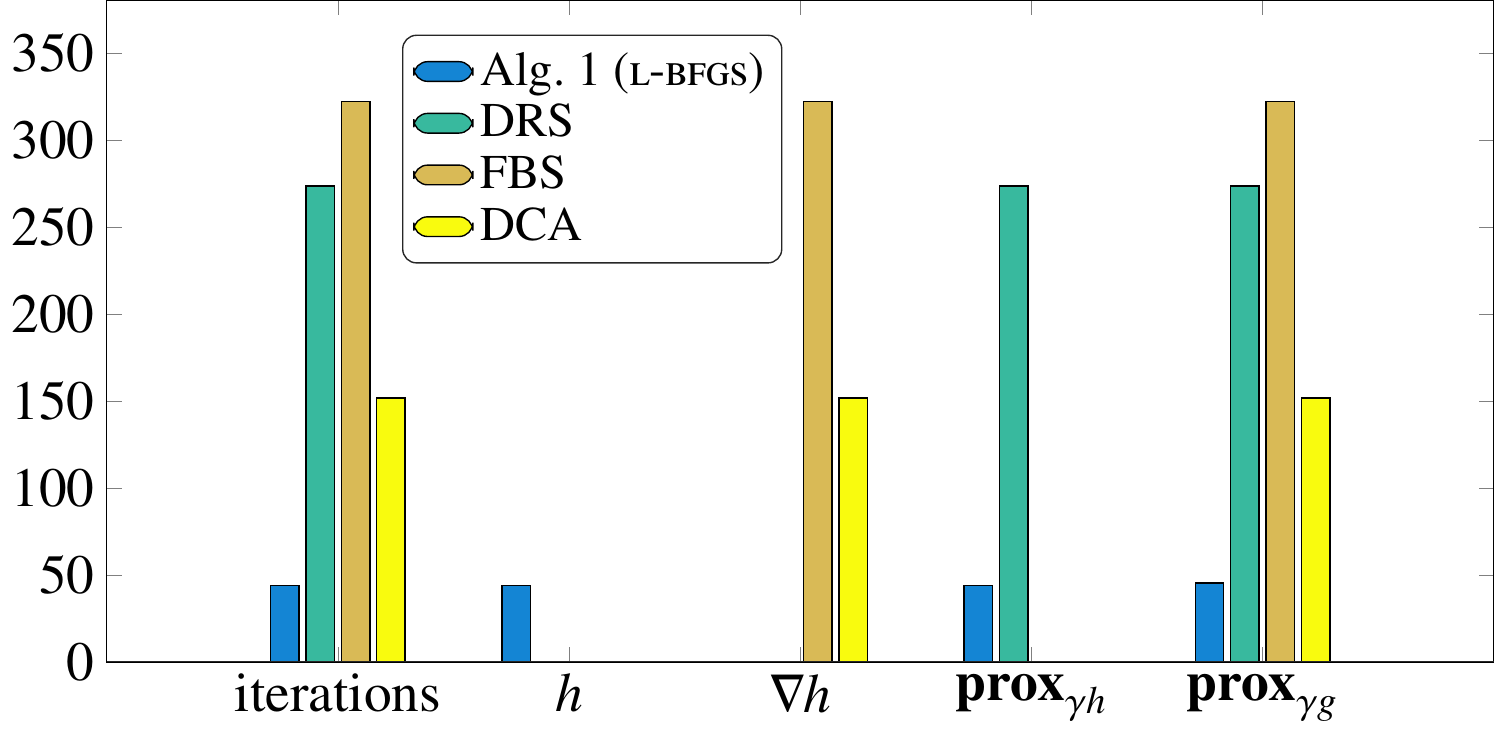}
		\caption{%
			Iteration comparison for random instances of \eqref{eq:PCA}.%
		}%
		\label{fig:iterations}
	\end{figure}

	To \eqref{eq:PCA} we applied FBS, DRS, DCA and \cref{alg:P} (gradient descent on the DCE) with L-BFGS steps and Wolfe backtracking.
	Sparse random matrices \(A\in\R^{20n\times n}\) with 10\% nonzeros were generated for 11 values of \(n\) on a linear scale between 100 and 1000, with a sufficiently small \(\kappa\) \cite[\S2.1]{journee2010generalized}.
	The mean number of iterations required by the solvers over these instances is reported in the first column of \cref{fig:iterations}.
	A stepsize \(\gamma = 0.9\lambda_\textrm{max}^{-1} (\Sigma) \) was selected for \cref{alg:P} and FBS, and \(\gamma = 0.45\lambda_\textrm{max}^{-1} (\Sigma) \) for DRS consistently with the nonconvex analysis in \cite{themelis2020douglas}.
	Stepsize tuning might lead to a better performance of these algorithms but was not considered here. 
	The termination criterion \(\|\prox_{\gamma\h}(s) - \prox_{\gamma\g}(s)\| \leq 10^{-6}\) was used for all solvers.
	Plain \cref{alg:P} (without L-BFGS) always exceeded 1000 iterations.
	
	\Cref{fig:iterations} also lists the complexity in terms of function calls.
	Evaluating \(\h\) and \(\nabla\h\) requires a matrix-vector product, which is \(O(n^2)\) operations.
	By factorizing \(\I +\gamma\Sigma\) once offline, each backsolve to compute \(\prox_{\gamma\h}\) also requires \(O(n^2)\) operations.
	Finally, \(\prox_{\gamma\g}\) requires \(2n\) comparisons and a norm-operation, and is clearly the least expensive operation.
	
	DCA and FBS need one \(\nabla\h\) and one \(\prox_{\gamma\g}\) (or similar) operation, and DRS one \(\prox_{-\gamma\h}\) (work equivalent to \(\prox_{\gamma\h}\)) and one \(\prox_{\gamma\g}\) operation per iteration.
	\cref{alg:P} requires one \(\prox_{\gamma\h}\) and one \(\prox_{\gamma\g}\) operation per iteration, and L-BFGS needs additionally one call to \(\h\), \(\prox_{\gamma\h}\) and \(\prox_{\gamma\g}\) per trial stepsize in the linesearch.
	However, as \(\h\) and \(\prox_{\gamma\h}\) involve linear operations for this particular problem, only one evaluation is required during the whole linesearch.
	Furthermore, in practice, it was observed that a stepsize of 1 was almost always accepted.
	From \cref{fig:iterations} it follows, therefore, that \cref{alg:P} with L-BFGS requires less work to converge than the other methods, disregarding the one time factorization cost not present in FBS and DCA.

	\section{Conclusions}\label{sec:Conclusions}
	By reshaping nonsmooth DC problems into the minimization of the smooth DC envelope function (DCE), a gradient method yields a new algorithm for DC programming.
	The algorithm is of splitting type, involving (subgradient-free, proximal) operations on each component which, additionally, can be carried out in parallel at each iteration.
	The smooth reinterpretation naturally leads to the possibility of Newton-type acceleration techniques which can significantly affect the convergence speed.
	The DCE has also a theoretical appeal in its deep kinship with the forward-backward envelope, as it is shown to be a reparametrization with more favorable reguarity properties.
	We believe that this connection may be a valuable tool for relaxing assumptions in FBE-based algorithms, which is planned for future work.


	\bibliographystyle{plain}
	\bibliography{TeX/Bibliography_init.bib}

\end{document}